\documentclass[aps,prl,twocolumn, superscriptaddress,showpacs,secnumarabic]{revtex4-1}
\usepackage{graphicx,xr,relsize,docs}
\usepackage{psfrag}
\usepackage{subfigure}
\usepackage{verbatim}
\usepackage{color}
\usepackage{amsmath,amssymb}
\newcommand{\bea}{\begin{eqnarray}}
\newcommand{\eea}{\end{eqnarray}}
\newcommand{\Xt}{\tilde{X}}
\newcommand{\Yt}{\tilde{Y}}

\newcommand{\Xb}{\mathbf{X}}
\newenvironment{claim}[1]{\par\noindent{\bf Claim:}\space#1}{}
\newenvironment{claimproof}[1]{\par\noindent{\bf Proof:}\space#1}{{\leavevmode\unskip\penalty9999 \hbox{}\nobreak\hfill\quad\hbox{$\blacksquare$}}}
\newtheorem{thm}{Theorem}[section]
\begin{document}

\title{Improving extreme value statistics}

\author{Ashivni Shekhawat}
\affiliation{Miller Institute for Basic Research in Science, UC Berkeley, Berkeley, CA USA}
\affiliation{Materials Sciences Division, Lawrence Berkeley National Laboratory, Berkeley, CA, USA}

\begin{abstract}
The rate of uniform convergence in extreme value statistics is non-universal and can be arbitrarily slow. 
Further, the relative error can be unbounded in the tail of the approximation, leading to difficulty in 
extrapolating the extreme value fit beyond the available data. We show that by using simple 
nonlinear transformations the extreme value approximation can be rendered rapidly convergent 
in the bulk, and asymptotic in the tail, thus fixing both issues. 
The transformations are often parameterized by just one parameter which can be estimated numerically.
The classical 
extreme value method is shown to be a special case of the proposed method.
We demonstrate that vastly improved results 
can be obtained with almost no extra cost. 
\end{abstract}



\maketitle
Extreme value statistics  provides a universal
\footnote{In a broad sense the term `universal' refers to any behavior 
that is largely independent of details, and is shared across a broad class. For example, 
the mean of iid random variables has a Gaussian distribution under very general conditions. 
This behavior is largely independent of the details of the distribution of the iid variables, and 
thus is universal. In similar spirit, the generalized extreme value distribution is a universal 
description of large statistical fluctuations.}
statistical description of 
rare events. Such events dictate the fate of a vast variety 
of phenomena spanning science, engineering, and humanities. 
Examples include stock 
price fluctuations~\cite{longin1996,broussard1998}, hydrology and one-hundred-year floods~\cite{katz2002,saralees2007,khaliq2006}, 
catastrophic fracture~\cite{weibull1939,harlow1978,duxbury1987,shekhawat2012}, 
climatology~\cite{buishand1989,smith1989},
risk management~\cite{resnick1999,mcNeil2000}, large insurance claims~\cite{beirlant1992,nader1997},
novelty detection~\cite{roberts2000} and so on.
The mathematical model for a rare event is that of a  
large statistical fluctuation. Let $\Xt$ be a random variable
with a cumulative distribution function (cdf) $F(\cdot)$, thus, $P(\Xt< x) = F(x)$.
Let $(\Xt_1,\ldots,\Xt_n)$ be a sample of $n$ iid random variables drawn 
from the distribution $F(\cdot)$. 
For example, these might represent the magnitude of the annual floods in years $(1,\ldots,n)$.
If one is interested in the largest flood, one might ask for its
cdf, i.e, $P(X_n< x)$, where $X_n=\mathrm{max}(\Xt_1,\ldots,\Xt_n)$, which obviously  
is given by $F(x)^n$. The pivotal result
of extreme value theory is that the distribution of maximum (or minimum) of 
iid random variables converges to a universal form under suitable 
linear rescaling. For a distribution $F(\cdot)$ if 
there exists a suitable sequence of constants, $a_n \in \Re,\ b_n \geq 0$
such that $F(a_n x + b_n)^n \to G_\gamma(x)$, then $G_\gamma(x)$ is 
of the form $\exp\{ - (1 + \gamma x)^{-1/\gamma} \}$, where $\gamma \in \Re$,
and $G_0(x) \equiv \exp\{-e^{-x}\}$~\cite{gnedenko1943,haan1970,resnick1987extreme}. 
A distribution function $F(\cdot)$ that satisfies the above for a given $\gamma$ is said 
to be in the domain of attraction of $G_\gamma(\cdot)$, or $F \in D(G_\gamma)$. 
The cases $\gamma\ =,\ >,\ <\ 0$ correspond to the
Gumbel (or type I), the Frechet (or type II), and the Weibull (or type III) distributions, respectively.
The conditions for $F(\cdot)$ to be in the domain of attraction of $G_\gamma(\cdot)$ are well 
established and fairly mild, see Refs.~\cite{resnick1987extreme,haan1996,kotz2000} for details.
Since the restrictions on $F(\cdot)$ are mild, this result is comparable to the central limit theorem  in its generality. 
However, the central limit theorem is a stronger result since the Berry-Essene theorem bounds
the rate of uniform convergence to the central limit under very general conditions (existence of
first three moments).
There is no analogous result in the theory of extremes.
\par
The success of extreme value theory is due to its simplicity and generality. Only three parameters, $\gamma,\ a_n,\ b_n$,
need to be fitted to data. Unfortunately, this stark simplicity is not carried over to the study of 
quality of approximation and rate of convergence. In classical extreme value theory, the rate of convergence
can vary widely, and needs to be evaluated on a case-by-case basis 
(see chapter 2 of ~\cite{resnick1987extreme} and Refs.~\cite{gomes1984,cohen1982,gyorgyi2010,anderson1978,anderson1984,haan1996}). 
It is the goal of this paper
to make the convergence properties more universal, at the cost of introducing slight complexity in the approximation
itself. We first discuss the issues associated with rate of convergence in the classical setting and then 
present the proposed formulation.
\par
There are two measures of quality of convergence that are considered widely. 
Firstly, one aims to bound the absolute maximum error of approximation,
$d_n = \sup_{x} | F(a_n x + b_n)^n - G_\gamma(x) |$.
The analogous bounds for the central limit theorem are provided by Berry-Essene type results.
Results of comparable generality are not available in the theory of extremes. 
Instead, the bound and its asymptotic behavior for large $n$ are to be evaluated on a case-by-case basis, 
and depend on the details of the tail of $F(\cdot)$ (see Ref.~\cite{resnick1987extreme} section 2.4
and supplemental sections~\ref{sup-sec:rc}, \ref{sup-sec:ex} for details).
The error of approximation can also be quantified via Edgeworth type expansions, which assert 
\begin{equation}
\lim_{n\to \infty} \frac{F(a_n x + b_n)^n -G_\gamma(x)}{W(n)} = G'_\gamma (x) \hat{H}_\gamma[G_\gamma(x)],
\label{eq:Edgeworth}
\end{equation}
uniformly in $x$, where the exact form of the function $\hat{H}(\cdot)$ is somewhat complicated (see Ref.~\cite{haan1996}). 
For the Edgeworth type expansions, the rate of convergence is governed by the $F(\cdot)$-dependent function $W(\cdot)$ (Eq.~\ref{eq:W}).
In either case, the decay of  $d_n$ and $W(\cdot)$ can be arbitrarily slow (or arbitrarily fast) depending on the 
tail properties of $F(\cdot)$. 
For example, $W(n) = 0$ if $F(x) = \exp( -e^{-x})$ (the Gumbel distribution), 
$W(n) \sim -1/2n$ if $F(x) = 1 - e^{-x}$ (the standard exponential distribution), and $W(n) \sim -1/2\log(n/\sqrt{2\pi})$
if $F(x) = \int_{-\infty}^{x} e^{-t^2/2}/\sqrt{2\pi} dt$ (the standard normal distribution); similar trends work for $d_n$ 
(supplemental section~\ref{sup-sec:ex}). 
Thus, the rate of convergence can range from infinitely fast to logarithmically slow (or worse). 
The logarithmic rate of convergence is obviously a cause of concern in practice. 
Refs.~\cite{gomes1984,cohen1982} show that the convergence can sometimes be improved by considering penultimate approximations, but the rate still remains 
logarithmic in several cases of interest.
In this paper we will show that the rate of convergence can be improved considerably ($1/n$ as opposed to $1/\log n$ )
in a robust and feasible manner. 
\par
A second measure of convergence has to do with the fact that the maximum error 
is not always the best measure of how close 
$F(a_nx + b_n)^n$ is to $G_\gamma(x)$ in the upper (or lower) tails.
The relative error in the tails is important for cases where one is interested
in the probability of large exceedances. 
In such cases, the quality of the upper tail of the approximation is 
measured by the ratio~\cite{anderson1978}
\begin{equation}
L(x) \equiv ({1-F(a_n x+ b_n)^n})/({1-G_\gamma(x)}).
\label{eq:LargeDev}
\end{equation}
Ideally $L(x)$ should stay close to 1. However, practically it 
can differ significantly from its ideal value of 1 for $x$ close to 
$x_+ \equiv \sup_x \{x: F(x) < 1\}$
at fixed $n$.
This behavior is characterized by studying the speed at which (for a given $n$) $x_n$ can be 
let to go to $x_+$, such that $L(x) \to 1$ uniformly for $x < x_n$~\cite{anderson1978}.
Here we take a more simple minded approach and study $\lim_{x\to\infty} L(x)$.
As before, there is a whole range of possible behavior. $L(x)$ approaches its 
ideal value of 1 for the exponential distribution, while it decays to 0 rather quickly
for the normal distribution. This behavior can lead to 
particularly severe errors and uncertainty when the fit to the extreme value 
approximation need to be extrapolated beyond the available data. 
This is typical 
of a large number of applications, such as prediction of large floods,
insurance claims or wild fires. Indeed, practitioners routinely predict the probability of
1000 year floods based on less than a century worth of good data! We will show how this difficulty
can be alleviated in our setup. 
\par
It is sometimes indicated in 
the literature that the slow rate of convergence is limited to the functions in domain of convergence 
of the Gumbel (type I) distribution, i.e.,~the cases where $F \in D(G_0)$. This is incorrect. Ref.~\cite{haan1996} shows that 
$F\in G_\gamma$ if the derivative of the function $j(x) \equiv F^{-1}(e^{-1/x})$ is regularly varying \footnote{A function
$f(x)$ is said to the regularly varying with index $\gamma$ at $\infty$ if $\lim_{t\to\infty} f(tx)/f(t) = x^{\gamma}$. It is called 
slowly varying if $\lim_{t\to\infty} f(tx)/f(t) = 1$.} with index $\gamma -1$, i.e.,
$\lim_{t\to\infty} j'(tx)/j'(t) = x^{\gamma -1}$. Since $j'(x)$ is regularly varying with index $\gamma -1$, it 
admits the representation $j'(x) = x^{\gamma -1}U(x)$, where $U(\cdot)$ is slowly varying, i.e., $\lim_{t\to\infty} U(tx)/U(t) = 1$. 
Note that the domain of convergence is solely controlled by the regularly varying part, $x^{\gamma -1}$, 
in the decomposition of $j'(x)$ and is independent of the slowly varying part, $U(x)$.
The rate of convergence is related to $j(\cdot)$ by
\begin{equation}
W(n) = nj''(n)/j'(n) - \gamma + 1 = nU'(n)/U(n), 
\label{eq:W}
\end{equation}
and is thus controlled solely by $U(\cdot)$, independent of $\gamma$. This shows that the convergence can be 
arbitrary irrespective of the domain of attraction. However, it is true that out of the most 
commonly used distributions, those belonging to $D(G_0)$ are more prone to such issues. We shall restrict 
our discussion to such case from here onwards. However, our method is equally applicable to other 
cases.
\par
The crux of our suggested methodology can be demonstrated by a simple example.
Consider an (admittedly contrived) industrial process that grinds out metallic disks
whose area, $A$, is distributed exponentially, so that $P(A<a) = 1-e^{-a}$.
Two analysts are given 100 boxes, each containing 10 such disks. They are  
asked to approximate the probability distribution of the radii of largest disk in each of the boxes.
The first analyst (say, Bob) simply measures the radius of the 
largest disk in each box, and fits these 100 observation to an extreme value 
form, perhaps using maximum likelihood estimation (MLE). 
The second analyst (say, Alice) decides to take a different route.
She measures the areas of the largest disk in each box, and fits this data to an extreme value
distribution instead, then she can predict the required probability by using a simple transformation.
Both Bob and Alice report their findings, and the employers who know the exact distribution fire Bob.
What went so wrong for Bob? 
\par
\begin{figure}[tbp]
\begin{center}
\includegraphics[width=0.4\textwidth,angle=0]{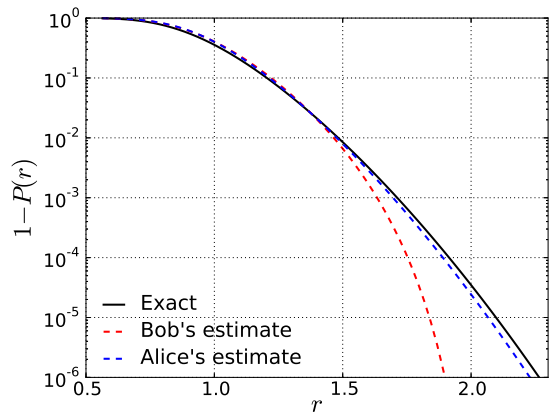}
\end{center}
\caption{Comparison of the exact and estimated probability of the radius of the largest of 10 disks being greater than $r$.
Bob's estimate, based on the extreme value theory, largely underestimate the probability of observing disks with large radii,
while Alice's estimate, based on treating area as the primary random variable, and transforming back the results to get probability 
of radii, works much better.} 
\label{fig:Radius}
\end{figure}
A great deal of insight about convergence issues in extreme value statistics can be gained 
by analyzing this example. Let $P(A<a),\ P(R<r)$ be the probabilities that the area, radius of a disk are lesser than $a$, $r$,
respectively. Clearly, $P(A<a) = 1 - e^{-a}$ and $P(R<r) = 1 - e^{-\pi r^2}$. It is also clear that the tail of radius distribution, $P(R<r)$, 
decays faster than exponentially for large $r$, thus any extreme value distribution, $G_\gamma(\cdot)$, will not be able to model it accurately. 
On the other hand, the tail of the area distribution, $P(A<a)$, decays exponentially, and can be modeled accurately by $G_0(\cdot)$. Thus, there 
is an inherent advantage to working with $A$ as the random variable being fit to extreme value distributions, even though $A$ and $R$ are simply related 
by $A = \pi R^2$. After a fit has been obtained for $A$, the probability for $R$ can be obtained easily by transforming back via $P(R<r) = P(A<\pi r^2)$.
Figure~\ref{fig:Radius} shows a comparison of Bob's and Alice's estimates and the exact result. 
Since there were 100 boxes, the empirical data was available at a probability level of $1-P = 10^{-2}$. Up to this level both estimates agree reasonably with the 
exact result. However, at $r = 1.91$, the exact result is $1-P(r) = 10^{-4}$, Bob's estimate is $1-P_B(r) = 4.5 \times 10^{-7}$, and Alice's estimate is $1-P_A(r) = 7.3 \times 10^{-5}$.
Thus, Alice's estimate is off by about 25\%, while Bob's is off by more than two orders of magnitude. Formally, one can show that  
$W(n) \sim 1/2n$, $\sim 1/\log n$, for Alice's and Bob's estimates, respectively. 
\par
The insight gained from the above example can be formalized. 
The idea is that it can be advantageous to work with a suitably transformed variable,
instead of the raw data itself.
The extreme value estimate for the raw data, $P(X_n< x) \approx G_{\gamma_x}( (x - b_{nx})/a_{nx})$ 
is susceptible to all the convergence issues discussed previously. 
\begin{claim}
There exists a monotonic $n$-independent transformation $\Yt = \hat{T}(\Xt)$ and constant $\gamma$ such that the extreme value approximation is exact, 
i.e.~$P(X_n < x)  = P(Y_n < \hat{T}(x)) = G_{\gamma}( (\hat{T}(x) - b_{ny})/a_{ny})$ for suitable
$n$-dependent constants $b_{ny},\ a_{ny}$, where $Y_n = \mathrm{max}(\Yt_1,\ldots,\Yt_n)$.
\end{claim}
\begin{claimproof}
$\hat{T}(x) = G_0^{-1}(F(x))$, $\gamma = 0$, $a_{ny} = 1$, $b_{ny} = \log n$ are suitable, as can be checked by direct substitution. 
However, this choice is not unique.
\end{claimproof}
\par
\noindent
Thus, working with a suitably transformed variable $\Yt = \hat{T}(\Xt)$ completely suppresses the systematic 
errors of the extreme value approximation in the sense of Eqs.~\ref{eq:Edgeworth},~\ref{eq:LargeDev}.
However, there is a slight problem with this scheme: it demands that to construct $\hat{T}(\cdot)$
we know $F(\cdot)$, which if we knew, we could calculate $F(\cdot)^n$ exactly without this elaborate scheme anyway!
This problem is made tractable by the following results.
\begin{claim}
Let $F(\cdot)$ have unbounded support (the case of bounded support is similar). Let  
$F(x) \sim 1 - \sum_{i=0}^{\infty}f_i(x)$ be an asymptotic expansion for large $x$, where the gauge functions
$f_i(\cdot)$ are monotonic. Then the variable $\Yt = T(\Xt) = -\log f_0(\Xt)$
is asymptotically exponentially distributed, and for the Edgeworth expansion 
corresponding to the variable $\Yt$ (Eq.~\ref{eq:Edgeworth}) the rate of convergence,
and the quality of the upper tail are characterized by 
\begin{equation}
W(n) \sim 1/2n,
\label{eq:convCl}
\end{equation}
\begin{equation}
(1-P(Y_n< y + \log n))/(1-G_0(x)) \to 1.
\label{eq:LargeDevCl}
\end{equation}
\end{claim}
\begin{claimproof}
Since $f_0(\cdot)$ is monotonic, $T^{-1}(x) = f_0^{-1}(e^{-x})$. Now, 
$P(\Yt < y) = P (\Xt < T^{-1}(y)) = F( f_0^{-1}(e^{-y}))\sim 1 - e^{-y}.$
Thus, the distribution of $\Yt$ is asymptotically exponential. 
Eqs.~\ref{eq:convCl},~\ref{eq:LargeDevCl} hold due to properties of 
the standard exponential distribution (see supplemental section~\ref{sup-sec:pr} for detailed proof).
\end{claimproof}
Thus, instead of knowing $F(\cdot)$ it is sufficient 
to estimate $f_0(\cdot)$. 
Since $P(Y_n<y) = P(X_n < T^{-1}(y))$, we get the following convergence assurances based 
on Eqs.~\ref{eq:convCl},~\ref{eq:LargeDevCl}
\begin{equation}
\lim_{n\to\infty} \frac{ F(T^{-1}(x + \log n))^n - G_0(x) }{1/2n} = G'_0(x) \hat{H}_0[G_0(x)],
\label{eq:assurConv}
\end{equation}
\begin{equation}
({1-F(T^{-1}(x+ \log n))^n})/({1-G_0(x)}) \to 1,
\label{eq:assurLargeDev}
\end{equation}
where $\hat{H}_0[G_0(x)] = e^{-x} + x -1$ (supplemental section~\ref{sup-sec:pr}). 
We have taken the norming constants $a_{ny},\ b_{ny} = 1,\ \log n$, 
as these are the theoretical asymptotic values for the exponential distribution. 
In practice they must be treated as free parameters to be fit. 
\par
The proposed method, which we call the T-method (`T' for transformation), is now clear.
Let us say that we can parameterize the transformation $T(x) = -\log f_0(x)$ by a parameter
$\beta$, then we have a parameter vector $\theta = (\beta, a_n, b_n)$, and a model
$G_0((T(\Xb|\beta) - b_n)/a_n)$. Given data vector $\Xb = (X_{n_1},X_{n_2},\ldots,X_{n_m})$, 
the parameter vector $\theta$ can be estimated 
via the maximum likelihood method by maximizing the following likelihood function
\begin{equation}
\mathcal{L}(\theta|\Xb) = \sum_{i=1}^m \big( \log G'_0( (T(\Xb_i) - b_n)/a_n) + \log T'(\Xb_i)/a_n\big).
\label{eq:MLE}
\end{equation}
\par
As a final step, the transformation $T(\cdot)$ needs to be parameterized by 
the parameter $\beta$ in a principled way. As mentioned 
previously, 
we restrict our discussion to distributions in the domain of $G_0(\cdot)$, 
and a transformation of the form discussed next will be useful
only if for the raw data we get $\gamma$ close to 0. The required transformation can be worked out 
easily for several common distributions with $F \in D(G_0)$. For example (supplemental section~\ref{sup-sec:tr}), for the normal distribution
we get $f_0(x) = e^{-x^2/2}/\sqrt{2\pi}x$, $T(x) = -\log f_0(x) \sim x^2$ (strictly speaking $T(x)\sim x^2/2 + \log(x/\sqrt{2\pi}) $, however
multiples can be absorbed into $a_{n}$, constants into $b_{n}$, and we ignore the asymptotically smaller $\log x$), 
for Rayleigh type distributions $F(x) = 1 - e^{-x^\alpha}$, $T(x) = x^\alpha$,
for lognormal distribution $T(x) \sim (\log x)^2$ etc. The heuristic is that 
if a semilog plot of $1-\hat{C}_n(x)$, where $\hat{C}_n(x)$ is the empirical cdf of the data, is a 
straight line, then the underlying distribution $1-F(x)$ is exponential, and no correction 
is needed. If the plot curves downwards, then $1-F(x)$ decays super-exponentially, and $T(x) = x^\beta$
with $\beta > 1$. If the semilog plot curves upwards, while a loglog plot curves downwards, then 
the decay is super-polynomial, but sub-exponential, thus $T(x)$ is of the form $x^\beta$ with $0<\beta < 1$ 
or of the form $\log(x)^\beta$. 
If the loglog plot is roughly straight, then it is likely that $ F\notin D(G_0)$, 
and either a correction is not needed or it is more subtle, and will be discussed in a later paper. Once a 
form is chosen, the parameter set $\beta,\ a_{n},\ b_{n}$ can be obtained by using MLE estimator
suggested in Eq.~\ref{eq:MLE} or another estimator in the usual manner. 
Note that $\gamma = 0$ is held fixed, so there are still only three free parameters in the model.
The classical extreme value fit is a special case of the T-method 
obtained when $T(\cdot)$ is taken to be the identity, i.e.~$T(x) = x$.
\begin{figure}[tbp]
\begin{center}
\includegraphics[width=0.4\textwidth,angle=0]{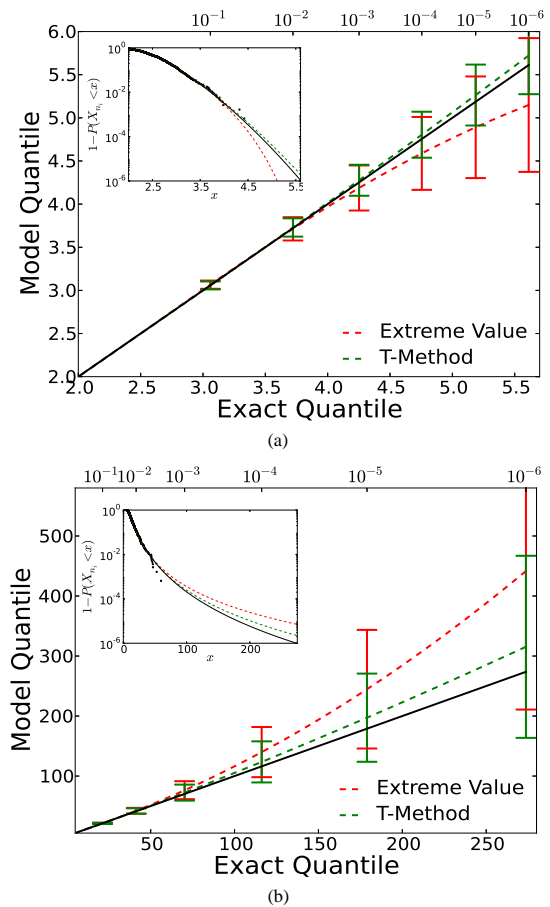}
\end{center}
\caption{Comparison of the classical extreme value approximation with the suggested transformation based method 
for data taken from (a) the normal distribution, and (b) the lognormal distribution. 
The main graph shows the traditional QQ plot with the upper $x$-axis showing
the exceedance, $1-P(X_{n_i})$ corresponding to the quantile on the main $x$-axis. The solid black line is a guide to 
the eye and shows the ideal result. The dashed lines show the model quantiles
averaged over 1000 monte carlo runs; while the errorbars show the 2-standard deviation
range. In each monte carlo run the model fit to a sample of size $m=1000$, 
and the fit is extrapolated to a probability level of $1-P(X_{n_i}<x) = 10^{-6}$.
The insets show the upper tail of the estimation, $1-P(X_{n_i} < x)$, on a semilog plot
for a typical monte carlo run; the empirical data is shown in the black dots.
It is clear that the transformation based method yields better predictions
and less variance
even when extrapolated well beyond the range of the available data. 
}\label{fig:Trans}
\end{figure}
\par
We test the proposed method on data generated from normal and lognormal distributions. For the case of 
the normal (lognormal) distribution, 
we generate a random sample $\mathbf{X} = (X_{n_1},\ldots,X_{n_m})$, where $m=1000$. Each $X_{n_i} = \mathrm{max}(\Xt_1,\ldots,\Xt_n)$
where $n=100$, and $\Xt_i$ are iid random variables drawn from the normal (lognormal) distribution. 
We estimate the parameter vector $\mathbf{\theta} = (\beta, a_{n}, b_{n})$ by using MLE (Eq.~\ref{eq:MLE})
with $T(x|\beta) = x^\beta$ for the normal case, and $T(x|\beta) = (\log x)^\beta$ for 
the lognormal case. Figure~\ref{fig:Trans} shows
a favorable comparison of the results obtained by the T-method with the classical extreme value approximation.
We have also tested the method on other distributions, including Rayleigh type distributions, and the 
Leath-Duxbury distribution encountered in statistics of fracture~\cite{duxbury1987}.
It is clear that the suggested method out-performs the classical extreme value approximation with the same number 
of parameters. Finally, we tested the method on the exponential distribution, where the 
convergence to $G_0(\cdot)$ is rapid, and a correction is not needed per se. We found that 
the T-method increases the mean accuracy of predictions slightly,
while reducing the spread in the predictions significantly (supplemental section~\ref{sup-sec:ft}).
\par
In summary, we have suggested a simple method, which we call the T-method, to alleviate the problem of slow convergence 
of classical extreme value approximations. 
The method works by estimating simple nonlinear transformation that defines a new random variable that has 
better convergence properties in the extreme value sense. 
Some previous authors have studied rates of 
convergence in nonlinear scaling in 
extreme value statistics (see Refs.~\cite{pantcheva1985,barakat2010}). Their results are rather 
remarkable, however, their focus has been on 
studying $d_n$ or $W(n)$ for specific transformations (power transformation, for example) rather 
than constructing numerical methods of wide applicability. In this sense the proposed T-method
is complementary to their results.
The T-method was applied succesfully to distributions in the domain of 
attraction of the Gumbel (type I) distribution.
We hope that application of our method will lead to more reliable estimates of probabilities of 
extremes in a large number of applications.
\begin{acknowledgments}
The author acknowledges support from DOE-BES
DE-FG02-07ER46393 while he was at Cornell University, and from
the Miller Institute for Basic Research in Science at University of California, Berkeley.
The author would like to thank Prof.~James P.~Sethna for insightful discussions, and 
Prof.~Robert O.~Ritchie for hosting him at the Lawrence Berkeley National Laboratory.
\end{acknowledgments}
\bibliography{improvingEVS}
\section{Supplemental Material}
\subsection{Domain of Attraction}
\label{sup-sec:doa}
Let $F(\cdot)$ be cumulative distribution function (cdf). A theorem due to 
Gnedenko (and modified by others to get various characterizations) is stated below;
see Refs.~\cite{resnick1987extreme,kotz2000,coles2001} for proofs.
\begin{thm}
If there exists a sequence of normalizing constants $a_n > 0$ and $b_n \in \Re$, such that
\begin{equation}
F(a_n x + b_n)^n \to G_\gamma(x), 1+\gamma x > 0,
\label{sup-eq:EVS}
\end{equation}
weakly as $n \to \infty$, then $G_\gamma(x)$ is of the form
\begin{equation}
G_\gamma(x) = \exp\{ - (1 + \gamma x)^{-1/\gamma}\},\quad \gamma \in \Re.
\end{equation}
\end{thm}
In such a case we say that $F(\cdot)$ is in the domain of attraction of $G_\gamma(\cdot)$, or $F \in D(G_\gamma)$.
A characterization of the domain of attraction is as follows.
\begin{thm}
(See Ref.~\cite{haan1996} for details) Let 
\begin{equation}
j \equiv \left(\frac{1}{-\log F}\right)^{-1},
\end{equation}
where the $(\cdot)^{-1}$ denotes inverse, i.e., under suitable continuity conditions $j(x) \equiv F^{-1}(e^{-1/x})$.
Then Eq.~\ref{sup-eq:EVS} holds iff
\begin{equation}
\lim_{t\to\infty} \frac{j(tx)-j(t)}{tj'(t)} = \frac{x^\gamma - 1}{\gamma}.
\end{equation}
A sufficient condition for Eq.~\ref{sup-eq:EVS} to hold is
\begin{equation}
\lim_{t\to\infty} \frac{j'(tx)}{j'(t)} = x^{\gamma - 1}
\end{equation}
\end{thm}

\subsection{Norming Constants}
\label{sup-sec:nc}
\begin{thm}
(See Ref.~\cite{haan1996} for details) 
The following constants are asymptotically optimal norming constants in Eq.~\ref{sup-eq:EVS}
\begin{equation}
a_n = n j'(n),\quad b_n = j(n).
\end{equation}
\end{thm}

\subsection{Rates of Uniform Convergence and Edgeworth Expansions}
\label{sup-sec:rc}
Let $F \in D(G_\gamma)$, and
\begin{align}
d_n &= \sup_n| F(a_n x + b_n)^n - G_\gamma(x)|,\\
x_0 & = \sup\left\{ x : F(x) < 1 \right\}, \\
h(x) &= -\log F(x) - \left\{ \frac{-F(x)F''(x)\log F(x)}{(F'(x))^2} + 1\right\}.
\end{align}
Let $g(x)$ be such that 
\begin{equation}
|h(x)| \leq g(x) \downarrow 0\ \mathrm{as}\ x \to x_0.
\end{equation}
\begin{thm}
(See Ref.~\cite{resnick1987extreme} section 2.4.2 for details) Let 
If $F \in D(G_0)$ then $d_n \leq O(g(b_n)).$
\end{thm}
\begin{thm}
\label{sup-th:EW}
(See Ref.~\cite{haan1996} for details) 
Let $F \in D(G_\gamma)$, then 
\begin{equation}
\lim_{n\to \infty} \frac{F(a_n x + b_n)^n -G_\gamma(x)}{W(n)} = G'_\gamma (x) \hat{H}_\gamma[G_\gamma(x)],
\label{sup-eq:Edgeworth}
\end{equation}
uniformly in $x$, where $\hat{H}_\gamma(G_\gamma(x)) = H_\gamma(-\log(-\log G_\gamma(x)))$ and
\begin{equation}
H_\gamma(x) = \left\{ \begin{array}{ll}
			\int_0^x e^{\gamma u} \int_0^u e^{\rho s} ds du & \mathrm{for}\ \gamma \geq 0\\
			-\int_x^\infty e^{\gamma u} \int_0^u e^{\rho s} ds du & \mathrm{for}\ \gamma < 0.
			\end{array}\right.,
\end{equation}
where $\rho$ is such that for $v(t) \equiv j(e^t)$, and  $A(e^t) \equiv v''(t)/v'(t) - \gamma$,
\begin{equation}
\lim_{t\to\infty} \frac{A(tx)}{A(t)} = x^\rho\ \mathrm{for}\ x > 0.
\end{equation}
The function $W(n)$ is given by
\begin{equation}
W(n) = nj''(n)/j'(n) -\gamma + 1.
\end{equation}
\end{thm}
\subsection{Examples}
\label{sup-sec:ex}
\subsubsection{Exponential Distribution}
\label{sup-sec:exED}
For the standard exponential distribution, $F(x) = 1 - e^{-x}$, $F^{-1}(x) = -\log(1-x)$. We get
\begin{equation}
j(x) = F^{-1}(e^{-1/x}) \sim \log x + \frac{1}{2x} + O(1/x^2). 
\end{equation}
Thus,
\begin{equation}
a_n = nj'(n) \sim 1,\quad b_n = j(n) \sim \log n.
\end{equation}
Further
\begin{align}
h(x) & = -\log F(x) - \left\{ \frac{-F(x)F''(x)\log F(x)}{(F'(x))^2} + 1\right\}\\
			& = -1 - e^x \log( 1- e^{-x})\\
			& = \frac{e^{-x}}{2} + \frac{e^{-2x}}{3} + \frac{e^{-3x}}{4} + \ldots\\
			& <  e^{-x}.
\end{align}
Thus we get $g(x) = e^{-x}$ and $g(b_n) = 1/n$, so that 
\begin{align}
d_n = \sup_n| F(a_n x + b_n)^n - G_\gamma(x)| \leq O(1/n)
\end{align}
Grinding through the calculations further gives
$A(t) \sim 1 /(2t - 1)$, thus, $\rho = -1$. Thus,
\begin{equation}
H_\gamma(x) = \int_0^x e^{\gamma u} \int_0^u e^{\rho s} ds du = e^{-x} + x - 1
\end{equation}
which yields 
\begin{equation}
\hat{H}_\gamma(x) = e^{-x} + x - 1.
\end{equation}
The function $W(n)$ is 
\begin{equation}
W(n) = nj''(n)/j'(n) -\gamma + 1 \sim 1/2n.
\end{equation}
Thus, the Edgeworth expansion becomes 
\begin{equation}
\lim_{n\to \infty} \frac{F(a_n x + b_n)^n -G_0(x)}{1/2n} = G'_0(x) (e^{-x} + x - 1)
\end{equation}

\subsubsection{Normal Distribution}
\label{sup-sec:exND}
\begin{equation}
F(x) = \frac{1}{\sqrt{2\pi}} \int_{-\infty}^x e^{-t^2/2} dt \sim 1 - \frac{e^{-x^2/2}}{x\sqrt{2\pi}} \left( 1 - \frac{1}{x^2}\right)
\end{equation}
Since $F(j(x)) = e^{-1/x}$, so for large $x$, $j(x)$ must be suitably large. However, $j(x)$ might not grow rapidly enough with 
$x$, so we keep a higher order term for $j(x)$ in the following expansion
\begin{align}
F(j(x)) & = e^{-1/x} \\
1 - \frac{e^{-j(x)^2/2}}{j(x)\sqrt{2\pi}} \left( 1 - \frac{1}{j(x)^2}\right) & \approx 1 - \frac{1}{x} \\
j^2/2 + \log j - \log (1 - \frac{1}{j^2}) \approx \log(x/\sqrt{2\pi}) \\
j^2/2 + \log j + 1/j^2 \approx \log(x/\sqrt{2\pi}) 
\end{align}
Ignoring the $1/j^2$ and solving gives
\begin{equation}
j(x) \sim \sqrt{2\log \tilde{x} -\log(2\log \tilde{x})},
\end{equation}
where $\tilde{x} = x/\sqrt{2\pi}$.
Grinding through the details, we get 
\begin{gather}
b(n) \sim  \sqrt{2\log \tilde{n} -\log(2\log \tilde{n})},\\
a(n) \sim 1/b(n),\ W(n) \sim -1/2\log(n).
\end{gather}
Further calculations show that $A(t) \sim -1/2\log t$, giving $\rho = 0$ and  
\begin{equation}
H_\gamma(x) = \int_0^x e^{\gamma u} \int_0^u e^{\rho s} ds du = x^2/2
\end{equation}
which yields 
\begin{equation}
\hat{H}_\gamma(x) = x^2/2.
\end{equation}
Thus, the Edgeworth expansion becomes 
\begin{equation}
\lim_{n\to \infty} \frac{F(a_n x + b_n)^n -G_0(x)}{-1/2\log n} = G'_0(x) \frac{x^2}{2}
\end{equation}

\subsubsection{LogNormal Distribution}
\label{sup-sec:exLND}
\begin{align}
F(x) & = \frac{1}{\sqrt{2\pi}} \int_{-\infty}^{\log x} e^{-t^2/2} dt\\
			& \sim 1 - \frac{e^{-(\log x)^2/2}}{\log x\sqrt{2\pi}} \left( 1 - \frac{1}{(\log x)^2}\right).
\end{align}
The analysis proceeds in a manner analogous to the last section. The first order results are 
\begin{equation}
a(n) \sim e^{D(n)}/D(n),\ b(n) \sim e^{D(n)},\ W(n) \sim -1/D(n),
\end{equation}
where $D(n) = (2\log(n/\sqrt{2\pi}))^{1/2}$.
\subsubsection{Rayleigh Distribution}
\label{sup-sec:exRD}
\begin{equation}
F(x) = 1 - e^{-x^\alpha},\quad F^{-1}(x) = (-\log(1-x))^{1/\alpha}
\end{equation}
Auxiliary function $j$
\begin{align}
j(x) & = F^{-1}(e^{-1/x}) \\
		 & \sim (\log x)^{1/\alpha}\left( 1 + \frac{1}{2\alpha x\log x}\right)
\end{align}
Thus
\begin{align}
b(n) & \sim (\log n)^{1/\alpha} \\
a(n) & \sim \frac{(\log n)^{1/\alpha -1}}{\alpha}\\
W(n) & \sim \frac{1-\alpha}{a\log n} + \frac{1}{2n}
\end{align}
\subsubsection{Gamma Distribution}
\label{sup-sec:exGD}

\begin{align}
F(x) & = \frac{1}{\Gamma(a)}\int_0^x t^{a-1} e^{-t}dt \\
		& \sim 1 - \frac{x^{a-1}e^{-x}}{\Gamma(a)}\left( 1 + \frac{a-1}{x} + \frac{(a-1)(a-2)}{x^2} + \ldots\right)
\end{align}
Thus
\begin{align}
j(x) & \sim \log (x/\Gamma(a)) + (a-1)\log(\log(x/\Gamma(a)))\\
b(n) & \sim \log (n/\Gamma(a)) + (a-1)\log(\log(x/\Gamma(a)))\\
a(n) & \sim 1 + (a-1)/\log(x/\Gamma(a))\\
W(n) & \sim -\frac{a-1}{(a-1 + \log(x/\Gamma(a)))\log(x/\Gamma(a))}
\end{align}

\subsection{Transformations}
\label{sup-sec:tr}
This section has the calculation for the asymptotic terms in the mapping $G^{-1}_0(F(x))$ for some $F(x)$.
\subsubsection{Normal}
\label{sup-sec:trND}
The normal cdf is,
\begin{equation}
F(x) = \frac{1}{\sqrt{2\pi}} \int_{-\infty}^x e^{-t^2/2} dt \sim 1 - \frac{e^{-x^2/2}}{x\sqrt{2\pi}} \left( 1 - \frac{1}{x^2}\right).
\end{equation}
while the inverse of $G_0(\cdot)$ is
\begin{align}
G^{-1}_0(x) & = -\log(-\log(x)) \\
G^{-1}_0(1-x) & \sim -\log(x) - \frac{x}{2} -\frac{5x^2}{24} - \ldots.
\end{align}
Thus, the transformation becomes
\begin{equation}
G^{-1}_0(F(x)) \sim \frac{x^2}{2} + \log\left( x\sqrt{2\pi} \right) + \mathcal{O}(1/x^2)
\end{equation}
\subsubsection{Lognormal }
\label{sup-sec:trLND}
The lognormal cdf is
\begin{equation}
F(x) \frac{1}{\sqrt{2\pi}} \int_{-\infty}^{\log x} e^{-t^2/2} dt.
\end{equation}
Thus, the transformation becomes
\begin{equation}
G^{-1}_0(F(x)) \sim \frac{(\log x)^2}{2} + \log\left( \log x\sqrt{2\pi} \right) + \mathcal{O}(1/\log x^2)
\end{equation}
\subsubsection{Rayleigh}
\label{sup-sec:trRD}
The Rayleigh cdf is
\begin{equation}
F(x) = 1 - e^{-x^\alpha}
\end{equation}
Thus the transformation becomes
\begin{equation}
G^{-1}_0(F(x)) \sim x^\alpha -\frac{e^{-x^\alpha}}{2} + \ldots
\end{equation}
\subsubsection{Gamma}
\label{sup-sec:trGD}
The Gamma cdf is 
\begin{align}
F(x) & = \frac{1}{\Gamma(a)}\int_0^x t^{a-1} e^{-t}dt \\
		& \sim 1 - \frac{x^{a-1}e^{-x}}{\Gamma(a)}\left( 1 + \frac{a-1}{x} + \frac{(a-1)(a-2)}{x^2} + \ldots\right)
\end{align}
Thus, the transformation becomes
\begin{equation}
G^{-1}_0(F(x)) \sim x + (1-a)\log x + \log \Gamma(a)
\end{equation}
\subsection{Tail Convergence}
\label{sup-sec:tc}
\subsubsection{Exponential Distribution}
\label{sup-sec:tcED}
The exponential cdf and norming constants are 
\begin{equation}
	F(x) = 1 - e^{-x},\ a_n = 1,\ b_n = \log n.
\end{equation}
Thus,
\begin{align}
L(x) & = \frac{1-F(a_n x + b_n)^n}{1-G_0(x)}\\
			& = \frac{1-(1 - e^{-x}/n)^n}{1-e^{-e^{-x}}}\\
			& = \frac{e^{-x} + o(e^{-x})} { e^{-x} + o(e^{-x})} \to 1,
\end{align}
where $o()$ is the `small-$o$' notation. Thus, the extreme value approximation for the exponential distribution is good in the upper tail 
of the distribution.

\subsubsection{Normal Distribution}
\label{sup-sec:tcND}
The normal cdf and norming constants are 
\begin{equation}
	F(x) = \frac{1}{\sqrt{2\pi}} \int_{-\infty}^x e^{-t^2/2} dt \sim 1 - \frac{e^{-x^2/2}}{x\sqrt{2\pi}} \left( 1 - \frac{1}{x^2}\right),
\end{equation}
\begin{equation}
b(n) \sim  \sqrt{2\log \tilde{n} -\log(2\log \tilde{n})},\ a(n) \sim 1/b(n),\ \tilde n = n/\sqrt{2\pi}.
\end{equation}

Thus, it is easy to see that 
\begin{align}
L(x) & = \frac{1-F(a_n x + b_n)^n}{1-G_0(x)} \to 0
\end{align}
Thus, the extreme value approximation for the normal distribution is bad in the upper tail.
\subsection{Proofs}
\label{sup-sec:pr}
\begin{claim}
Let $F(\cdot)$ have unbounded support (the case of bounded support is similar). Let  
$F(x) \sim 1 - \sum_{i=0}^{\infty}f_i(x)$ be an asymptotic expansion for large $x$, where the gauge functions
$f_i(\cdot)$ are monotonic. Then the variable $Y = T(X) = -\log f_0(x)$
is asymptotically exponentially distributed, i.e.,
\begin{equation}
P(Y < y) \sim 1 - e^{-y},
\end{equation}
and
\begin{equation}
\lim_{n\to\infty}  \frac{P(Y_n < y + \log n) - G_0(y)}{W(n)} = G'_0(y)\hat{H}_0[G_0(y)],
\label{sup-eq:convCl}
\end{equation}
where $ W(n) \sim 1/2n$,
and
\begin{equation}
(1-P(Y_n< y + \log n))/(1-G_0(x)) \to 1.
\label{sup-eq:LargeDevCl}
\end{equation}
\end{claim}
\begin{claimproof}
Since $f_0(\cdot)$ is monotonic, $T^{-1}(x) = f_0^{-1}(e^{-x})$. Now, 
\begin{align} 
K(y) \equiv P(Y < y) 	& = P (X < T^{-1}(y)) \\
											& = F( f_0^{-1}(e^{-y}))\\
											& \sim 1 - e^{-y}.
\end{align}
Thus, the distribution of $Y$ is asymptotically exponential. 
Consider the auxiliary function $j(n)$. By definition
\begin{equation}
K(j(n)) = e^{-1/n},
\end{equation}
Thus, $j(n) \to \infty$ as $ n \to \infty$, since it is easy to check that $K(y)$ is monotonic and $\lim_{y\to\infty} K(y) = 1$.
Thus, for large $n$ we can use the asymptotic expansion $K(j(n))\sim 1 - e^{-j} + o(e^{-j})$. The monotonity of 
the gauge functions $f_i(\cdot)$ ensure that there are no oscillatory terms in this asymptotic expansion, and thus 
it can be differentiated term by term. Thus, we have
\begin{equation}
1 - e^{-j(n)} + o(e^{-j(n)}) \sim e^{-1/n}.
\end{equation}
The above has the solution 
\begin{equation}
j(n) \sim -\log ( 1 - e^{-1/n}) + o(1/n),
\end{equation}
as can be verified by direction substitution. 
Thus, for the Edgeworth expansion (theorem~\ref{sup-th:EW}), we get 
\begin{equation}
W(n) = nj''(n)/j'(n) - \gamma + 1 = nj''(n)/j'(n) + 1 \sim 1/2n.
\end{equation}
Further following theorem~\ref{sup-th:EW}, we get $A(t) \sim 1/2t$, thus giving $\rho = -1$. Since $\gamma = 0$, 
we get $\hat{H}_0[G_0(x)] = e^{-x} + x + 1$ from theorem~\ref{sup-th:EW}. The Edgeworth expansion is thus established.
Finally, for the tail approximation
\begin{align}
L(y) 	& = \frac{1-P(Y_n< y + \log n)}{1-G_0(y)}\\
			& = \frac{1-K(y + \log n)^n}{1-G_0(y)}\\
			& = \frac{1-(1 - e^{-y}/n + o(e^{-y}/n) )^n}{1-e^{-e^{-y}}}\\
			& = \frac{e^{-y} + o(e^{-y}) }{e^{-y} + o(e^{-y})} \to 1,
\end{align}
and the proof is complete.
\end{claimproof}
\subsection{Fits to Exponential Data}
\label{sup-sec:ft}
\begin{figure}[tbp]
\begin{center}
\includegraphics[width=0.4\textwidth,angle=0]{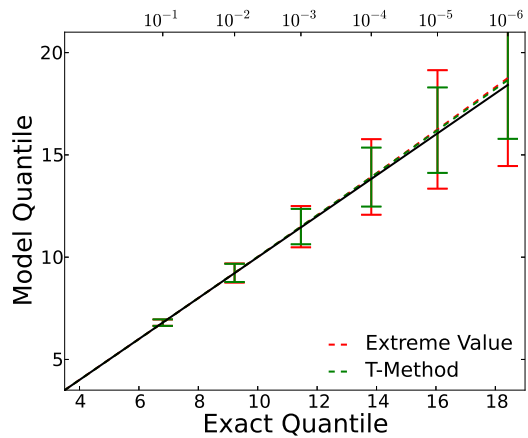}
\end{center}
\caption{Comparison of the classical extreme value approximation with the suggested transformation based method 
for data taken from the standard exponential distribution. 
The solid black line is a guide to 
the eye and shows the ideal result. The dashed lines show the model quantiles
averaged over 1000 monte carlo runs; while the errorbars show the 2-standard deviation
range. In each monte carlo run the model fit to a sample of size 1000, 
and the fit is extrapolated to a probability level of $1-P(X_n<x) = 10^{-6}$.
The insets show the upper tail of the estimation, $1-P(X_n < x)$, on a semilog plot
for a typical monte carlo run; the empirical data is shown in the black dots.
It is clear that the transformation based method yields better predictions
and less variance
even when extrapolated well beyond the range of the available data. 
}\label{sup-fig:SupTrans}
\end{figure}

Here we apply the T-method to the standard exponential distribution, $F(z) = 1 - e^{-z}$.
Since the convergence of the exponential distribution to the 
extreme value form $G_0(\cdot)$ is rapid in the bulk as well as in the tail, 
the application of the T-method is not necessary to get a good fit to the 
data, or a good result from the extrapolation of the fit. However, we  
apply the method to test if its application in such a case result
in predictions that are any worse (or better) than the standard 
extreme value statistics. In particular, we consider the 
distribution of the variable $X$ = $\mathrm{max}(Z_1,\ldots,Z_m)$,
where $Z_i$ are exponential iid random variables, $P(Z_i < z) = F(z) = 1 - e^{-z}$.
We take $m = 100$, and take a sample $\mathbf{X} = (X_1,\ldots,X_n)$ of size
$n = 1000$. The classical extreme value model is parameterized
by the parameter vector $\theta_c = (\gamma,a_n,b_n)$, and leads
to the following log-likelihood function
\begin{equation}
\mathcal{L}_c(\theta_c|\mathbf{X}) = \mathlarger{\Sigma}_i \big( \log G'_\gamma((X_i - b_n)/a_n) + \log (\gamma/a_n) \big).
\end{equation}
While, the T-method with the transformation $Y = T(X|\beta) = X^\beta$ is parameterized by
the parameter vector $\theta_t = (\beta,a_n, b_n)$ and leads to the following log-likelihood
function
\begin{multline}
\mathcal{L}_t(\theta_t|\mathbf{X}) = \mathlarger{\Sigma}_i \big(\log G'_0((T(X_i|\beta) - b_n)/a_n) \\
			+ \log T'(X_i|\beta) - \log a_n\big).
\end{multline}
We do monte carlo simulations by generating 1000 samples $\mathbf{X}$ and fitting the models. 
Figure~\ref{sup-fig:SupTrans} shows the mean and the 2-standard deviation bounds for the QQ plots of the fits. 
It is clear that the mean prediction from T-method is slightly better than the classical 
extreme value fit, while the T-method results in smaller variance in the predictions.
\end{document}